\documentclass[12pt]{amsart}
\usepackage{amsmath,amscd,amssymb,amsfonts}
\newcommand{\nin}{\par\noindent}
\newcommand{\bs}{\par\bigskip}
\newcommand{\ms}{\par\medskip}
\newcommand{\sk}{\par\smallskip}
\newcommand{\mb}{\raise.27ex\h{${\scriptscriptstyle\bullet}$}}
\newcommand{\ssb}{\raise.32pt\h{${\scriptscriptstyle\bullet}$}}
\newcommand{\sssb}{\raise-2.7pt\h{${\scriptscriptstyle\bullet}$}}
\newcommand{\ssc}{\,\raise.2ex\hbox{${\scriptstyle\circ}$}\,}
\newcommand{\sotim}{\hbox{\1$\otimes$}\1}
\newcommand{\motim}{\hbox{$\bigotimes$}}
\newcommand{\mopl}{\hbox{$\bigoplus$}}
\newcommand{\msum}{\hbox{$\sum$}}
\newcommand{\mprod}{\hbox{$\prod$}}
\newcommand{\mcup}{\hbox{$\bigcup$}}
\newcommand{\mcap}{\hbox{$\bigcap$}}
\newcommand{\mwedge}{\hbox{$\bigwedge$}}
\newcommand{\mf}{{\mathfrak m}}
\newcommand{\la}{\langle}
\newcommand{\ra}{\rangle}
\newcommand{\bl}{\bigl}
\newcommand{\br}{\bigr}
\newcommand{\af}{{\mathfrak a}}
\newcommand{\AAA}{{\bf A}}
\newcommand{\A}{{\mathcal A}}
\newcommand{\AX}{{\mathcal A}_X}
\newcommand{\At}{\widetilde{\mathcal A}}
\newcommand{\B}{{\mathcal B}}
\newcommand{\C}{{\mathcal C}}
\newcommand{\CC}{{\bf C}}
\newcommand{\DD}{{\bf D}}
\newcommand{\D}{{\mathcal D}}
\newcommand{\DX}{{\mathcal D}_X}
\newcommand{\E}{{\mathcal E}}
\newcommand{\F}{{\mathcal F}}
\newcommand{\G}{{\mathcal G}}
\newcommand{\Gs}{{\mathcal G}_{\ssb}}
\newcommand{\Gt}{\widetilde{\mathcal G}}
\newcommand{\h}{\hbox}
\newcommand{\Hc}{{\mathcal H}}
\newcommand{\Hct}{\widetilde{H}}
\newcommand{\I}{{\mathcal I}}
\newcommand{\J}{{\mathcal J}}
\newcommand{\K}{{\mathcal K}}
\newcommand{\Ks}{{\mathcal K}_{\ssb}}
\newcommand{\Lc}{{\mathcal L}}
\newcommand{\Ls}{{\mathcal L}_{\ssb}}
\newcommand{\Lt}{\widetilde{\mathcal L}}
\newcommand{\M}{{\mathcal M}}
\newcommand{\Mr}{{}^{r\!}{\mathcal M}}
\newcommand{\Nc}{{\mathcal N}}
\newcommand{\Ncb}{\,\overline{\!{\mathcal N}}}
\newcommand{\N}{{\bf N}}
\newcommand{\Oc}{{\mathcal O}}
\newcommand{\OP}{{\mathcal O}_{\mathcal P}}
\newcommand{\OX}{{\mathcal O}_X}
\newcommand{\pb}{\bar{\pi}}
\newcommand{\PP}{{\mathcal P}}
\newcommand{\PPP}{{\bf P}}
\newcommand{\Q}{{\bf Q}}
\newcommand{\q}{\quad}
\newcommand{\R}{{\bf R}}
\newcommand{\Rb}{\,\overline{\!R}}
\newcommand{\rd}{\partial}
\newcommand{\Sb}{{}\,\overline{\!S}{}}
\newcommand{\X}{{\mathcal X}}
\newcommand{\Y}{{\mathcal Y}}
\newcommand{\Ybs}{\overline{Y}_{\!s}}
\newcommand{\Z}{{\bf Z}}
\newcommand{\coh}{{\rm coh}}
\newcommand{\Ker}{{\rm Ker}}
\newcommand{\DR}{{\rm DR}}
\newcommand{\Ext}{{\rm Ext}}
\newcommand{\cExt}{{\mathcal E}xt}
\newcommand{\Hcom}{{\mathcal H}om}

\newcommand{\PProj}{{\mathcal P}roj}
\newcommand{\Spec}{{\mathcal S}pec}
\newcommand{\Sing}{{\rm Sing}\,}
\newcommand{\Gr}{\hbox{\rm Gr}}
\newcommand{\1}{{\hskip1pt}}
\newcommand{\simto}{\buildrel\sim\over\longrightarrow}

\newcommand{\into}{\hookrightarrow}
\newcommand{\too}{\longrightarrow}
\newcommand{\Cech}{{\v C}ech}
\newcommand{\ges}{\geqslant}
\newcommand{\les}{\leqslant}
\begin{document} \title{Graded duality for filtered $D$-modules}
\author{Morihiko Saito}
\address{RIMS Kyoto University, Kyoto 606-8502 Japan}
\email{msaito@kurims.kyoto-u.ac.jp}
\author{Christian Schnell}
\address{Department of Mathematics, Stony Brook University, Stony Brook, NY 11794, USA}
\email{cschnell@math.sunysb.edu}
\begin{abstract}
For a coherent filtered D-module we show that the dual of each graded piece over the structure sheaf is isomorphic to a certain graded piece of the ring-theoretic local cohomology complex of the graded quotient of the dual of the filtered D-module along the zero-section of the cotangent bundle. This follows from a similar assertion for coherent graded modules over a polynomial algebra over the structure sheaf. We also prove that the local cohomology sheaves can be calculated by using the higher direct images of the twists of the associated sheaf complex on the projective cotangent bundle. These are closely related to local duality, essentially due to Grothendieck.
\end{abstract}
\maketitle
\centerline{\bf Introduction}
\bs\nin
Let $(\M,F)$ be a filtered holonomic left $\DX$-module on a smooth complex algebraic variety $X$ of dimension $n$.
If $\Gr^F_{\ssb}\M:=\mopl_p\,\Gr^F_p\M$ is Cohen-Macaulay over $\AX:=\Gr^F_{\ssb}\DX$ (for instance, if $(\M,F)$ underlies a mixed Hodge module, see \cite[Lemma~5.1.13]{mhp}), then its $\D$-dual $(\M',F):=\DD(\M,F)$ is a filtered $\DX$-module, see (2.2) below for $\DD(\M,F)$.
However, the relation between the graded quotients $\Gr^F_{\ssb}\M$ and $\Gr^F_{\ssb}\M'$ is not obvious, and there were no results about it in the literature before \cite{Sch1} as far as we know.
\sk
Assume first that the restriction of $(\M,F)$ over an open subvariety $U$ of $X$ is a filtered vector bundle.
In case $(\M,F)$ underlies a mixed Hodge module, it underlies a variation of mixed Hodge structure over $U$ where $F_p=F^{-p}$.
Let $\DD_{\Oc}(\M,F)|_U$ denote the $\Oc$-dual of $(\M,F)|_U$ viewed as a filtered vector bundle.
It is known that there is a canonical isomorphism
$$\DD(\M,F)|_U=\DD_{\Oc}(\M,F[-n])|_U,$$
where $(F[m])_p=F_{p-m}$ for $p,m\in\Z$.
For instance, if $\M=\OX$ with $\Gr^F_p\OX=0\,\,(p\ne 0)$ so that $\DD_{\Oc}(\OX,F)=(\OX,F)$, then $(\OX,F)$ has weight $n$, and hence
$$\DD(\OX,F)=(\OX,F)(n)\,(:=(\OX,F[n])).$$
\sk
By the above isomorphism we have a canonical isomorphism of graded $\A_U$-modules
$$\Gr^F_{\ssb}\M'|_U=\Hcom_{\Oc_U}(\Gr^F_{n-\ssb}\,\M|_U,\Oc_U),$$
where the action of the odd degree part of $\A_U$ on the right-hand side is twisted by a minus sign.
However, it is not clear how this is extended over $X$.
In this paper we show the following theorem as a special case of Theorem~2 below.
Set $\I:=\mopl_{k>0}\,\Gr^F_k\DX$.
This is the graded ideal sheaf of $\AX$ corresponding to the zero-section $X$ of the cotangent bundle $T^*X$.
\ms\nin
{\bf Theorem~1.} {\it With the above notation and assumptions, the above isomorphism is extended to a canonical isomorphism of graded $\AX$-modules
$$\Gamma_{\I}(\Gr^F_{\ssb}\M')\sssb=\Hcom_{\OX}(\Gr^F_{n-\ssb}\,\M,\OX),$$
where the action of $\AX$ on the right-hand side is induced by the one on $\Gr^F_{n-\ssb}\,\M$ and is twisted by a minus sign for the odd degree part of $\AX$.}
\ms
The left-hand side of the isomorphism is the ring-theoretic local cohomology which is defined to be the $\I$-primary torsion part of $\Gr^F_{\ssb}\M'$, i.e., the union of the $\I^k$-torsion part for $k\gg 0$, see \cite{Gro1}, \cite{Gro2}.
(We avoid the notation $\Gamma_{\I}(\Gr^F_{\ssb}\M')$ since $\Gamma_{\I}(\Gr^F_p\M')$ does not make sense.) If we take local coordinates $x_1,\dots,x_n$ of $X$ and set $\xi_i:=\Gr^F_1\rd/\rd x_i\in\Gr^F_1\DX$, then
$$\Gamma_{\I}(\Gr^F_{\ssb}\M')\sssb=\mcup_{k\gg 0}\bl(\mcap_{i=1}^n\, \Ker(\xi^k_i:\Gr^F_{\ssb}\M'\to\Gr^F_{\ssb+k}\M')\br).$$
\sk
Let $D^b_{\coh}F(\DX)$ be the derived category of bounded complexes of filtered left $\DX$-modules with coherent cohomology sheaves.
Let $D^bG(\AX)$ be the derived category of bounded complexes of graded $\AX$-modules.
\ms\nin
{\bf Theorem~2.} {\it For $(\M,F)\in D^b_{\coh}F(\DX)$, let $(\M',F):=\DD(\M,F)\in D^b_{\coh}F(\DX)$, the dual filtered complex.
Then there is $\R\Gamma_{\I}(\Gr^F_{\ssb}\M')\sssb$ in $D^bG(\AX)$ together with a canonical isomorphism
$$\R\Gamma_{\I}(\Gr^F_{\ssb}\M')\sssb= \R\Hcom_{\OX}(\Gr^F_{n-\ssb}\,\M,\OX),$$
where the action of $\AX$ on the right-hand side is induced by the one on $\Gr^F_{n-\ssb}\,\M$ and is twisted by a minus sign for the odd degree part of $\AX$.}
\ms
If $(\M,F)$ is a Cohen-Macaulay holonomic filtered $\DX$-module, then the local cohomology sheaves $\Hc^i_{\I}(\Gr^F_{\ssb}\M')\sssb$ can be calculated by using the higher direct images of the twists of the associated coherent sheaf on the projective cotangent bundle together with the long exact sequence in Proposition~(3.9) below.
In case the support of $\M$ is not $X$ nor a union of smooth subvarieties, it is quite nontrivial to calculate both sides of the formula in Theorem~2 even in a simple case where the support is a quadratic surface with an ordinary double point, see Example~(4.3) below.
\sk
Theorem~2 follows from the following similar assertion for graded $\AX$-modules.
Let $D^b_{\coh}G(\AX)$ be the full subcategory of $D^bG(\AX)$ consisting of complexes with coherent cohomology.
\ms\nin
{\bf Theorem~3.} {\it For $\Gs\in D^b_{\coh}G(\AX)$, there is $\R\Hcom_{\AX}(\Gs,\AX[n])_{\ssb}$ in $D^b_{\coh}G(\AX)$, which is denoted by $\Gs'$.
Moreover, there is $\R\Gamma_{\I}\1(\Gs')\sssb$ in $D^bG(\AX)$ together with a canonical isomorphism
$$\R\Gamma_{\I}\1(\Gs')\sssb=\R\Hcom_{\OX}(\G_{-n-\ssb},\omega_X),$$
where the action of $\AX$ on the right-hand side is induced by the one on $\G_{-n-\ssb}$ without sign.}
\ms
The difference between $\Gr^F_{n-\ssb}$ in Theorem~2 and $\G_{-n-\ssb}$ in Theorem~3 comes from the transformation between left and right $\D$-modules in (2.1.3) since Theorem~3 implies the assertion for right $\D$-modules.
Theorem~3 is a special case of Theorem~(1.5) below which holds for any vector bundle (instead of the cotangent bundle) over any locally Noetherian scheme or analytic space.
We have $\E=\Omega_X^1$ if $V$ is the cotangent bundle on a smooth variety $X$ in (1.1).
\sk
The statement of the main theorem in \cite{Sch1} was given by using a combination of the assertions of Theorem~2 and Proposition~(3.9) below.
This is useful for explicit calculations, see Examples~(4.1) and (4.3) below.
A statement equivalent to Theorem~3 was noted without proof in Lemma~4.2 as a reformulation of Proposition~2.5 in loc.~cit.\ which was stated in terms of a spectral sequence.
The assertion was essentially proved there for coherent filtered $\DX$-modules on smooth quasi-projective varieties.
So it is enough to construct a canonical morphism to get the assertion for any bounded complexes of filtered $\DX$-modules with coherent cohomology sheaves.
Here it is rather nontrivial to keep track of the grading for the quasi-coherent sheaves on the vector bundle $V$ corresponding to quasi-coherent graded $\AX$-modules, especially when we take the localization and also the local cohomology sheaf along the zero-section of the vector bundle $V$ in order to get a distinguished triangle in the derived category, see Remark~(3.11)(iii) below.
This is related to Theorem~(3.8) and Proposition~(3.9) below (and also to the proof of Theorem in \cite[Section 1.2]{Sch1}).
\sk
Theorem~3 and Proposition~(3.9) may have been known to A.~Grothendieck (see, for instance, \cite[Proposition~2.1.5]{Gro1} and also the intended table of contents in \cite[p.~81]{Gro1}).
In the non-graded case, assertions similar to Theorem~3 have been studied by many people including J.P.C. Greenlees, J.P. May, J. Lipman and others (see \cite{GM}, \cite{AJL}, etc.) where one has to use completion and the arguments are more complicated.
Theorem~3 seems to be also related to Martineau-Harvey duality (\cite{Ma}, \cite{Ha}) for the analytic local cohomology of the structure sheaf along one point of a complex manifold, see also \cite{Gro2}, \cite{RD}, \cite{Na}, etc.
\sk
The first-named author is partially supported by Kakenhi 24540039.
The second-named author is partially supported by NSF grant DMS-1331641.
\sk
In Section 1 we prove a generalization of Theorem~3 after explaining some basics about the ring-theoretic local cohomology.
In Section 2 we apply this to filtered $\D$-modules and prove Theorem~2.
In Section 3 we study the graded localization along the zero-section to get a distinguished triangle in the derived category of graded modules.
In Section 4 we calculate some explicit examples to understand Theorem~2 and Proposition~(3.9) more concretely.
In Section 5 we explain an application of Corollary~(3.10) in a more general situation than in \cite[Section 1.3]{Sch1}.
\bs\bs
\centerline{\bf 1.~Graded duality for graded modules}
\bs\nin
In this section we prove a generalization of Theorem~3 after explaining some basics about the ring-theoretic local cohomology.
\ms\nin
{\bf 1.1.~Graded algebra associated to a vector bundle.} Let $\pi:V\to X$ be a vector bundle of rank $n>0$, where $X$ is a locally Noetherian scheme or an analytic space.
Let $\E$ be the corresponding locally free $\OX$-module.
Set
$$\E^{\vee}:=\Hcom_{\OX}(\E,\OX),\q \A_{X,i}:={\rm Sym}^i_{\OX}\E^{\vee},\q \AX=\A_{X,\ssb}:=\mopl_{i\ges 0}\,\A_{X,i},$$
where ${\rm Sym}^i_{\OX}\E^{\vee}$ denotes the symmetric product, which is the maximal quotient of $\motim^i\E^{\vee}$ on which the symmetric group acts trivially.
Then $\AX$ is a graded algebra over $\OX$.
This is locally a polynomial algebra, and we have
$$V=\Spec_X\AX\,\,\,\h{or}\,\,\,\,{{\Spec}an}_X\AX.$$
Let $D^bG(\AX)$ be the derived category of bounded complexes of graded $\AX$-modules, and $D^b_{\coh}G(\AX)$ the full subcategory consisting of complexes with coherent cohomology.
We have $D^+G(\AX)$, $D^-G(\AX)$ by replacing ``bounded'' with ``bounded below'' or ``bounded above'', and similarly for $D^{\pm}_{\coh}G(\AX)$.
\sk
The notation in this section is compatible with the situation in the introduction by setting $V=T^*X$ and $\E=\Omega_X^1$.
\ms\nin
{\bf 1.2.~Graded local cohomology.} With the notation and the assumption of (1.1), set
$$\I=\I_{\ssb}:=\mopl_{i>0}\,\A_{X,i}\subset\AX=\A_{X,\ssb}.$$
For a graded $\AX$-module $\Gs$, the graded local cohomology sheaf along the graded ideal sheaf $\I=\I_{\ssb}$ is defined by
$$\Gamma_{\I}\,(\Gs)_{\ssb}:=\mcup_{k\gg 0}\,\Gs^{\,\I^k}\q\h{with}\q \Gs^{\,\I^k}:=\Hcom_{\AX}(\AX/\I^k,\Gs)_{\ssb}\subset\Gs,
\leqno(1.2.1)$$
where $\Gs^{\,\I^k}$ consists of local sections of $\Gs$ annihilated by $\I^k$.
\sk
Let $\Gs\to\Ks$ be an injective resolution, see Remark~(1.6) below.
Then the derived local cohomology sheaf complex $\R\Gamma_{\I}\1(\Gs)_{\ssb}\in D^+G(\AX)$ for $\Gs\in D^+G(\AX)$ is defined by
$$\R\Gamma_{\I}\1(\Gs)_{\ssb}:=\mcup_{k\gg 0}\,\Ks^{\,\I^k}\q \h{with}\q\Ks^{\,\I^k}:=\Hcom_{\AX}(\AX/\I^k,\Ks)_{\ssb}\subset\Ks,
\leqno(1.2.2)$$
Here $\I^k$ in (1.2.1--2) can be replaced by a decreasing sequence of coherent ideal sheaves $\J_k\,\,(k\gg 0)$, if there are sequences $a(k)$, $b(k)$ such that
$$\I^{a(k)}\subset\J_k,\q\J_{b(k)}\subset\I^k\,\,\,(k\gg 0).$$
\sk
Since $\Ks$ is an injective complex, we have a canonical isomorphism
$$\Ks^{\,\I^k}=\R\Hcom_{\AX}(A_X/\I^k,\Gs)_{\ssb}\q\h{in}\q D^+G(\AX),
\leqno(1.2.3)$$
i.e., the left-hand side is a canonical representative of the right-hand side.
This is needed to define the inductive limit (using the union).
\ms
For an integer $p$, we will denote by $(p)$ the shift of grading so that
$$\G(p)_{\ssb}=\G_{\ssb+p}.
\leqno(1.2.4)$$
This is not compatible with the shift of an increasing filtration $F$ via the graded quotient since the sign differs.
\ms\nin
{\bf 1.3.~Lemma.} {\it With the above notation there is a canonical isomorphism}
$$\R\Hcom_{\AX}(\AX/\I,\AX[n])_{\ssb}=\mwedge^n\E\,(n)\q\h{in}\q D^bG(\AX).
\leqno(1.3.1)$$
\ms\nin
{\it Proof.} Note first that $\mwedge^n\E\,(n)$ has degree $0$ as a complex and pure degree $-n$ as a graded module by the above definition (1.2.4).
We have a graded free resolution of $\AX/\I$ by the complex whose $i$-th component is
$$\AX\sotim_{\OX}\mwedge^{-i}\E^{\vee}(i)\q(i\in[-n,0]),
\leqno(1.3.2)$$
where the differential is naturally defined since $\A_{X,1}=\E^{\vee}$.
(Here $\mwedge^{-i}\E^{\vee}$ has pure degree $0$, and not $-i$, as a graded module.) Indeed, this is the Koszul complex if we locally trivialize $V$ so that $\AX=\OX[\xi_1,\dots,\xi_n]$ locally.
Then the left-hand side of (1.3.1) is represented by the complex whose $i$-th component is
$$\AX\sotim_{\OX}\mwedge^{n+i}\E(n+i)\q(i\in[-n,0]).
\leqno(1.3.3)$$
This is a Koszul complex up to the twist by $\mwedge^n\E(n)$, since we have a canonical isomorphism
$$\mwedge^{n+i}\E(n+i)=\mwedge^n\E(n)\sotim_{\OX} \mwedge^{-i}\E^{\vee}(-i).$$
So the assertion follows.
\ms\nin
{\bf 1.4.~Lemma.} {\it With the notation of {\rm (1.1)}, assume $V$ is trivial so that $\E$ is free and $\AX=\OX[\xi_1,\dots,\xi_n]$.
Then there is an isomorphism in $D^bG(\AX)$
$$\R\Gamma_{\I}\bl(\OX[\xi_1,\dots,\xi_n]\br)[n]= \OX[\xi_1^{-1},\dots,\xi_n^{-1}]\1(\xi_1\cdots\xi_n)^{-1},
\leqno(1.4.1)$$
where the $\xi_i^{-1}$ have degree $-1$, and the right-hand side is identified with a quotient of $\OX[\xi_1,\xi_1^{-1},\dots,\xi_n,\xi_n^{-1}]$ to define the action of $\AX=\OX[\xi_1,\dots,\xi_n]$.}
\ms\nin
{\it Proof.} We may assume $X={\rm Spec}\,\Z$ or ${\rm Specan}\,\CC$ since the assertion is reduced to this case using base change by the canonical morphism from $X$ to ${\rm Spec}\,\Z$ or ${\rm Specan}\,\CC$ (depending on whether $X$ is algebraic or analytic).
So $\OX$ in (1.4.1) can be replaced by $\Z$ or $\CC$ which will be denoted by $R$.
Set
$$\J_k:=(\xi_1^k,\dots,\xi_n^k)\subset\A:=R[\xi_1,\dots,\xi_n] \subset\At:=R[\xi_1,\xi_1^{-1},\dots,\xi_n,\xi_n^{-1}].$$
\sk
Consider the shifted Koszul complex
$$\C_k^{\ssb}:=K^{\ssb}\bl(\A;\xi_1^k,\dots,\xi_n^k\br)[n],$$
which is associated to the multiplication by $\xi_j^k\,\,(j\in[1,n])$ on $\A$.
We have an identification
$$\C_k^{-i}=\mopl_{|J|=i}\,\A\,\xi_J^k\q\h{with}\q \xi_J^k:=\mprod_{j\in J}\xi_j^k,$$
where the $J$ are subsets of $\{1,\dots,n\}$, and the differential is given by the sum of natural inclusions $\A\,\xi_J^k\into\A\,\xi_{J'}^k\,\,(J\supset J')$ with appropriate sings.
\sk
Since $\{\xi_1^k,\dots,\xi_n^k\}$ is a regular sequence, we have a quasi-isomorphism
$$\C^{\ssb}_k\simto\A/\J_k,
\leqno(1.4.2)$$
By (1.2.3) together with the remark after (1.2.2) applied to $\J_k$, the left-hand side of (1.4.1) is represented by the inductive limit for $k\to\infty$ of
$$\C_k^{\prime\,\ssb}:=\xi^{-k}\C_k^{\ssb}\q\h{with}\q \xi:=\mprod_{i=1}^n\xi_i.
\leqno(1.4.3)$$
Indeed, (1.2.3) for $\Gs=\A$ with $\I^k$ replaced by $\J_k$ is represented by (1.4.3) which is the dual over $\A$ of $\C^{\ssb}_k$ in (1.4.2).
\sk
By (1.4.2) we have the quasi-isomorphism
$$\C^{\prime\,\ssb}_k\simto(\A/\J_k)\1\xi^{-k},
\leqno(1.4.4)$$
and the canonical surjection $\A/\J_{k'}\to\A/\J_k\,\,(k<k')$ induces the canonical inclusion
$$(\A/\J_k)\1\xi^{-k}\into (\A/\J_{k'})\1\xi^{-k'},
\leqno(1.4.5)$$
since the differential of $\C_k^{\ssb}$ is induced by the natural inclusions of free $\A$-submodules of $\widetilde{\A}$ up to appropriate signs, and this is the same for the dual $\C_k^{\prime\,\ssb}$.
So the assertion follows from (1.4.4--5) by taking the inductive limit.
\sk
The following may be known to A.~Grothendieck (see the intended table of contents in \cite[p.~81]{Gro1}).
\ms\nin
{\bf 1.5.~Theorem.} {\it With the notation and the assumption of {\rm (1.1)}, let $\Gs\in D^b_{\coh}G(\AX)$.
Then we have
$$\Gs':=\R\Hcom_{\AX}(\Gs,\AX[n])_{\ssb}\q\h{in}\q D^+_{\coh}G(\AX),$$
where $n={\rm rank}\,V$.
Moreover, there is a canonical isomorphism in $D^+G(\AX)$
$$\R\Gamma_{\I}\1(\Gs')\sssb\simto \R\Hcom_{\OX}(\G_{-n-\ssb},\mwedge^n\E),
\leqno(1.5.1)$$
where the right-hand side is defined by taking an injective resolution of the $\OX$-module $\mwedge^n\E$, and the action of $\AX$ is induced by that on $\G_{-n-\ssb}$ without sign.}
\ms\nin
{\it Proof.} There is an injective resolution (see Remark~(1.6) below):
$$\AX[n]\to\Ks.$$
So $\R\Hcom_{\AX}(\Gs,\AX[n])_{\ssb}$ is represented by $\Hcom_{\AX}(\Gs,\Ks)_{\ssb}$ in $D^+G(\AX)$, We then define as a complex of graded $\AX$-modules
$$\Gs':=\Hcom_{\AX}(\Gs,\Ks)_{\ssb}.$$
We may assume that $\Gs$ is a complex of flat $\AX$-modules by replacing it with a flat resolution.
Then $\Gs'$ is an injective complex.
Hence we have a canonical isomorphism in $D^+G(\AX)$
$$\R\Gamma_{\I}\1(\Gs')\sssb=\Gamma_{\I}\1(\Gs')\sssb.$$
\sk
We have a canonical morphism in $C^+G(\AX)$
$$\Gamma_{\I}\1(\Gs')\sssb \to\Hcom_{\AX}\bl(\Gs,\Gamma_{\I}\1(\Ks)_{\ssb}\br){}_{\ssb}.
\leqno(1.5.2)$$
By Lemma~(1.3) and (1.2.2--3) for $\Gs=\AX$, we get a canonical morphism of $\OX$-modules
$$\mwedge^n\E\cong\cExt^n(\AX/\I,\AX)_{-n}\to\Hc^n_{\I}\1(\AX)_{-n}\, \bl(\cong\Gamma_{\I}\1(\Ks)_{-n}\br),
\leqno(1.5.3)$$
and this is an isomorphism by the local calculation in Lemma~(1.4).
Using the restriction, we then get canonical morphisms of complexes of graded $\OX$-modules
$$\aligned\Hcom_{\AX}\bl(\Gs,\Gamma_{\I}\1(\Ks)_{\ssb}\br){}_{\ssb} &\to\Hcom_{\OX}\bl(\G_{-n-\ssb},\Gamma_{\I}\1(\Ks)_{-n}\br)\\ &\to\R\Hcom_{\OX}\bl(\G_{-n-\ssb},\mwedge^n\E\br),\endaligned
\leqno(1.5.4)$$
where the last morphism can be defined by taking an injective resolution of $\Gamma_{\I}\1(\Ks)_{-n}$.
The compatibility with the action of $\AX$ follows from the definition of the action of $\AX$ on the first term.
Composing (1.5.2) with (1.5.4), we thus get the canonical morphism (1.5.1).
By construction the obtained morphism is functorial for $\Gs$ contravariantly.
\sk
Once the canonical morphism is constructed, the assertion is local, and we may assume that $V$ is a trivial bundle, i.e.\ $\E$ is a free $\OX$-module, and $X$ is a Noetherian affine scheme or a closed analytic subspace in a polydisk (which is extendable to a closed analytic subspace in a slightly larger polydisk).
So $\AX$ is a polynomial algebra $\OX[\xi_1,\dots,\xi_n]$
\sk
Using the canonical morphism (1.5.1) constructed above together with the truncation $\tau_{\les p}$ (see \cite[1.4.6]{De}) on $\Gs$, the assertion is reduced to the case $\Gs$ is isomorphic to a coherent graded $\AX$-module up to a shift in $D^-G(\AX)$.
Replacing $\Gs$ with a resolution, we may then assume that $\Gs$ is bounded above and each component is isomorphic to a finite direct sum of $\AX(-p_i)$, where $(-p_i)$ is as in (1.2.4).
Using the truncation $\sigma_{\ges p}$ (see \cite[1.4.7]{De}), the assertion is further reduced to the case of $\AX(-p_i)$, where we can apply Lemma~(1.4).
So the assertion follows.
This finishes the proof of Theorem~(1.5).
\ms\nin
{\bf 1.6.~Remarks.} (i) We can show that the category of graded $\AX$-modules has enough injective objects by the same argument as in the usual case.
Indeed, using the sheaf of discontinuous sections at each degree, this can be reduced to the case of graded $A$-modules $G_{\ssb}$ with $A=\A_{X,x}$.
We have a surjection
$$\widetilde{G}'_{\ssb}\to G'_{\ssb}:={\rm Hom}_{\Z}(G_{-\ssb},\Q/\Z),$$
taking generators $g_{\alpha}$ of $G'_{\ssb}$ over $A$ where $\widetilde{G}'_{\ssb}$ is the direct sum of $\A(-p_{\alpha})$ for $p_{\alpha}\in\Z$ with $g_{\alpha}\in G'_{p_{\alpha}}$, see (1.2.4) for $(-p_{\alpha})$.
Then we get an injection
$$G_{\ssb}\to{\rm Hom}_{\Z}(\widetilde{G}'_{-\ssb},\Q/\Z),$$
and the target is an injective graded $A$-module.
\ms
(ii) The above argument implies that the category of quasi-coherent graded $\AX$-modules has enough injective objects in the case $X$ is affine by applying it to $\A=\Gamma(X,\AX)$.
In general we take an affine covering $\{U_j\}$, and use the adjunction between the direct image and the pull-back since quasi-coherent sheaves are stable by direct images.
\ms
(iii) For an injective graded $\AX$-module $\K_{\ssb}$, each $\K_p$ is an injective $\OX$-module.
Indeed, we have for any $\OX$-module $\F$
$$\Hcom_{\OX}(\F(p),\Ks)= \Hcom_{\AX}(\AX\sotim_{\OX}\F(p),\Ks),$$
where $(p)$ for $p\in\Z$ is as in (1.2.4).
\ms\nin
{\bf 1.7.~Remark.} For a complex manifold $X$ of dimension $n$, it is well known that there is a perfect pairing between the topological vector spaces
$$\Hc^n_{\{x\}}\OX\q\h{and}\q\Oc_{X,x}=\CC\{\!\{x_1,\dots,x_n\}\!\},$$
which is a special case of Martineau-Harvey duality (this is closely related to the theory of hyperfunctions, see \cite{Ma}, \cite{Ha}).
Analogues in algebraic geometry have been studied by many people, see, for instance, \cite{Gro1}, \cite{AJL}, \cite{GM}, \cite{Na}, \cite{RD} among others.
\bs\bs
\centerline{\bf 2.~Application to filtered $\D$-modules}
\bs\nin
In this section we apply Theorem~3 to filtered $\D$-modules and prove Theorem~2.
\ms\nin
{\bf 2.1.~Transformation between left and right $\D$-modules.} From now on, assume $X$ is a smooth complex algebraic variety or a complex manifold of dimension $n$.
The Hodge filtration $F$ on the right $\DX$-modules $\Omega_X^n$ and $\omega_X$ is given by
$$\Gr^F_p\Omega_X^n=0\,\,(p\ne -n),\q\Gr^F_p\omega_X=0\,\,(p\ne 0).
\leqno(2.1.1)$$
This implies
$$(\Omega_X^n,F)=(\omega_X,F)(-n):=(\omega_X,F[-n]).
\leqno(2.1.2)$$
Here $(F[m])_p=F_{p-m}$ in general.
Note that $(\Omega_X^n,F)$ and $(\omega_X,F)$ have weight $n$ and $-n$ respectively, see \cite{mhp}
\sk
We have the transformation from filtered left $\D$-modules to filtered right $\D$-modules defined by
$$(\M,F)\mapsto(\Mr,F):=(\Omega_X^n,F)\sotim_{\OX}(\M,F),
\leqno(2.1.3)$$
where the action of a vector field $v$ on $\Omega_X^n\sotim_{\OX}\M$ is given by
$$(\zeta\sotim m)v:=\zeta v\sotim m-\zeta\sotim vm\q\h{for}\,\,\, \zeta\in\Omega_X^n,\,m\in\M.
\leqno(2.1.4)$$
The shift of the filtration $F$ is in order to get the isomorphism of the associated filtered de Rham complexes
$$\DR_X(\M,F)=\DR_X(\Mr,F),$$
where the left-hand side is defined as usual up to a shift of complex, and the right-hand side is defined by using $\mwedge^{-\ssb}\Theta_X$, see loc.~cit.
To indicate left or right $\D$-modules, we will note $^l$ or $^r$ at the end of $D^bF(\DX)$ (for instance, $D^bF(\DX)^l$ for left $\D$-modules).
Then (2.1.3) induces equivalences of categories
$$D^bF(\DX)^l\simto D^bF(\DX)^r,\q D^b_{\coh}F(\DX)^l\simto D^b_{\coh}F(\DX)^r,\,\,\h{etc}.
\leqno(2.1.5)$$
\ms\nin
{\bf 2.2.~Dual of filtered $\D$-modules.} Let $\omega_X[n]\to\K_X$ be an injective resolution of a complex of right $\DX$-modules, or more generally, a quasi-isomorphism of complexes of right $\DX$-modules such that each component of $\K_X$ is injective over $\OX$.
(Here we use the flatness of $\D_X$ over $\Oc_X$.) It has the filtration $F$ defined by $\Gr^F_p\K_X=0$ for $p\ne 0$ as in (2.1.1).
For $(\M,F)\in D^b_{\coh}F(\DX)^r$, the dual filtered complex is defined by
$$\DD(\M,F)=\Hcom_{\DX}\bl(\DR_X^{-1}\DR_X(\M,F), \K_X\sotim_{\OX}(\DX,F)\br).
\leqno(2.2.1)$$
Here $\DR_X^{-1}\DR_X(\M,F)$ denotes the induced filtered $\DX$-module associated to the filtered de Rham complex $\DR_X(\M,F)$, and its $i$-th component is
$$\bl((\M,F[-i])\sotim_{\OX}\mwedge^{-i}\Theta_X\br) \sotim_{\OX}(\DX,F).$$
Taking local coordinates $x_1,\dots,x_n$, this complex is identified with the Koszul complex associated to the action of $\rd/\rd x_i$ on $\M\sotim_{\OX}\DX$ which is a right $\DX$-bi-module, see loc.~cit.
\ms
For $(\M,F)\in D^b_{\coh}F(\DX)^l$, its dual $(\M',F)=\DD(\M,F)\in D^b_{\coh}F(\DX)^l$ is defined by combining (2.2.1) with (2.1.5) so that
$$(\Omega_X^n,F)\sotim_{\OX}(\M',F)= \DD\bl((\Omega_X^n,F)\sotim_{\OX}(\M,F))\br)\q\h{in}\,\,\, D^b_{\coh}F(\DX)^r.
\leqno(2.2.2)$$
\ms\nin
{\bf 2.3.~Proposition.} {\it For $(\M,F)\in D^b_{\coh}F(\DX)^r$, let $(\M',F)$ be the dual filtered complex $\DD(\M,F)\in D^b_{\coh}F(\DX)^r$.
Then we have a canonical isomorphism in $D^bG(\AX)$
$$\Gr^F_{\ssb}\M'=\R\Hcom_{\AX}(\Gr^F_{\ssb}\M, \omega_X\sotim_{\OX}\AX[n])_{\ssb},
\leqno(2.3.1)$$
where the action of the odd degree part of $\AX$ on $\omega_X\sotim_{\OX}\AX$ is twisted by a minus sign.}
\ms\nin
{\it Proof.} We have to compare the right-hand side of (2.3.1) with the graded quotient of the right-hand side of (2.2.1).
The twist by the sign of the action of the odd degree part of $\AX$ on $\omega_X\sotim_{\OX}\AX$ comes from (2.1.4).
Concerning (2.2.1), there is a canonical filtered quasi-isomorphism
$$\DR_X^{-1}\DR_X(\M,F)\simto(\M,F).$$
Using this, we may replace the representative, and assume that each component of $(\M,F)$ is an induced filtered $\DX$-module, i.e., of the form
$$(\Lc,F)\sotim_{\OX}(\DX,F),$$
where $(\Lc,F)$ is a filtered $\OX$-module satisfying $F_p\Lc=0$ for $p\ll 0$.
\sk
For a filtered right $\DX$-module $(\Nc,F)$, we always assume $F_p\Nc=0$ for $p\ll 0$, and we have
$$\Hcom_{\DX}\bl((\Lc,F)\sotim_{\OX}(\DX,F),(\Nc,F)\br) =\Hcom_{\OX}\bl((\Lc,F),(\Nc,F)\br).
\leqno(2.3.2)$$
If the $\Gr^F_p\Nc$ are injective $\OX$-modules, we have moreover
$$\Gr^F_{\ssb}\Hcom_{\OX}\bl((\Lc,F),(\Nc,F)\br) =\Hcom_{\OX}\bl(\Gr^F_{\ssb}\Lc,\Gr^F_{\ssb}\Nc\br).
\leqno(2.3.3)$$
The above hypothesis is satisfied for $\K_X^i\sotim(\DX,F)$ since $\K_X^i$ is injective and the $\Gr^F_p\DX$ are locally free $\OX$-modules.
We have also
$$\Gr^F_{\ssb}\bl((\Lc,F)\sotim_{\OX}(\DX,F)\br)= \Gr^F_{\ssb}\Lc\sotim_{\OX}\Gr^F_{\ssb}\DX= \Gr^F_{\ssb}\Lc\sotim_{\OX}\AX.
\leqno(2.3.4)$$
For an injective resolution $\omega_X\sotim_{\OX}\AX[n]\to\Ks$, we have
$$\Hcom_{\AX}(\Gr^F_{\ssb}\Lc\sotim_{\OX}\AX,\Ks^i)\sssb =\Hcom_{\OX}(\Gr^F_{\ssb}\Lc,\Ks^i).
\leqno(2.3.5)$$
Comparing the right-hand side of (2.3.1) with the graded quotient of the right-hand side of (2.2.1) and using (2.3.2--5), the assertion (2.3.1) is then reduced to the assertion that for any injective graded $\AX$-module $\Ks$, each graded piece $\K_j$ is an injective $\OX$-module, see Remark~(1.6)(iii).
So Proposition~(2.3) follows.
\ms\nin
{\bf 2.4.~Theorem.} {\it For $(\M,F)\in D^b_{\coh}F(\DX)^r$, let $(\M',F):=\DD(\M,F)\in D^b_{\coh}F(\DX)^r$.
Then there is a canonical isomorphism in $D^bG(\AX)$
$$\R\Gamma_{\I}\1(\Gr^F_{\ssb}\,\M')\sssb= \R\Hcom_{\OX}(\Gr^F_{-n-\ssb}\,\M,\omega_X)\sotim_{\OX}\omega_X,$$
where the action of $\AX$ on the right-hand side is induced by the one on $\Gr^F_{-n-\ssb}\,\M$ and is twisted by a minus sign for the odd degree part of $\AX$.}
\ms\nin
{\it Proof.} Note first that $\sotim_{\OX}\omega_X$ on the right-hand side of the formula is needed as is seen in case $\M=\M'=\omega_X$.
This theorem follows from Theorem~3 and Proposition~(2.3).
Indeed, set
$$\Gs=\Gr^F_{\ssb}\,\M,\q \Gs':=\R\Hcom_{\AX}(\Gs,\AX[n])_{\ssb},$$
where the latter is the same as in Theorem~3.
By Proposition~(2.3) we have
$$\Gr^F_{\ssb}\M'=\Gs'\sotim_{\OX}\omega_X,$$
where the action of the odd degree part of $\AX$ on the right-hand side is twisted by a minus sign.
Since $\omega_X$ is an invertible sheaf (i.e., a locally free $\OX$-module of rank one), the assertion then follows from Theorem~3.
\ms\nin
{\bf 2.5.~Proof of Theorem~2.} Apply Theorem~(2.4) to
$$(\Mr,F):=(\Omega_X^n,F)\sotim_{\OX}(\M,F),\q (\Mr',F):=(\Omega_X^n,F)\sotim_{\OX}(\M',F).$$
Then Theorem~2 follows from this together with (2.2.2) and (2.1.2).
\bs\bs
\centerline{\bf 3.~Graded localization along the zero-section}
\bs\nin
In this section we study the graded localization along the zero-section to get a distinguished triangle in the derived category of graded modules.
We assume $X$ is a locally Noetherian scheme in this section.
\ms\nin
{\bf 3.1.~Graded localization.} In order to keep track of the grading of the localization along the zero-section of a vector bundle in \cite{Sch1}, we introduce the following ring-theoretic localization functor for graded $\AX$-modules $\Gs$ in the notation of (1.1--2):
$$\Gs(*\I)_{\ssb}:=\rlap{\raise-9pt\h{$\,\,\scriptstyle k$}}\rlap{\raise-5.5pt\h{$\,\rightarrow$}}{\rm lim}\,\Hcom_{\AX}(\I^k_{\ssb},\Gs)_{\ssb},
\leqno(3.1.1)$$
where we have in the notation of (1.2.4)
$$\Hcom_{\AX}(\I^k_{\ssb},\Gs)_p:=\Hcom_{\AX}(\I^k_{\ssb},\Gs(p)).$$
\sk
For $\Gs\in D^bG(\AX)$, its derived functor, denoted by $\R\Gamma(*\I)(\Gs)_{\ssb}$, is defined by taking an injective resolution $\Gs\to\Ks$ and setting
$$\R\Gamma(*\I)(\Gs)_{\ssb}:=\rlap{\raise-9pt\h{$\,\,\scriptstyle k$}}\rlap{\raise-5.5pt\h{$\,\rightarrow$}}{\rm lim}\,\Hcom_{\AX}(\I^k_{\ssb},\Ks)_{\ssb}.
\leqno(3.1.2)$$
Comparing this with the definition of $\R\Gamma_{\I}$ in (1.1), we get a canonical distinguished triangle
$$\R\Gamma_{\I}\1(\Gs)_{\ssb}\to\Gs\to\R\Gamma(*\I)(\Gs)_{\ssb} \buildrel{+1}\over\to\q\h{in}\,\,\,D^bG(\AX).
\leqno(3.1.3)$$
In \cite{Ka} an analogous construction in the analytic case was used.
Combining (3.1.3) with Theorem~(3.7) below, we get a distinguished triangle in Theorem~(3.8) which is somehow related to the proof of Theorem in \cite[Section 1.2]{Sch1}, see Remark~(3.11)(iii) below.
\ms\nin
{\bf 3.2.~Associated sheaf.} Set $\PP=\PProj_X\AX$.
Let $\OP(1)$ be the associated ample invertible sheaf on $\PP$, and set $\OP(p):=\OP(1)^{\otimes p}$ as usual.
Let $\pb:\PP\to X$ denote the natural projection.
We denote by $\Gt$ the quasi-coherent $\OP$-module associated to a quasi-coherent graded $\AX$-module $\Gs$.
Set $\Gt(p):=\Gt\sotim_{\OP}\OP(p)$.
This is identified with $\widetilde{\G(p)}$ where $(p)$ denote the shift of grading such that $\G(p)_i:=\G_{p+i}$, see (1.2.4).
Note that for $g\in\Gamma(U,\A_{X,i})$ with $U$ an affine open subvariety of $X$, the sections of $\Gt(p)$ over $\pb^{-1}(U)\setminus g^{-1}(0)$ is given by the degree $k$ part of the graded localization of $\Gamma(U,\Gs|_U)$ by $g$.
\sk
Let
$$R^j\pb_*\Gt(\mb):=\mopl_{p\in\Z}\,R^j\pb_*\Gt(p).$$
We have a canonical morphism of graded $\AX$-modules
$$\Gs(*\I)_{\ssb}\to\pb_*\Gt(\mb):=\mopl_{p\in\Z}\,\pb_*\Gt(p).
\leqno(3.2.1)$$
Indeed, a graded $\AX$-linear morphism $\I^k_{\ssb}\to\Gs(p)$ induces a morphism of quasi-coherent $\OP$-modules
$$\OP\to\Gt(p),$$
and this induces (3.2.1) by considering the image of the unit section $1$, since the inclusions $\I^k_{\ssb}\into\I^{k'}_{\ssb}$ for $k>k'$ induce the identity $\OP\to\OP$.
This is compatible with action of
$$\A_X=\mopl_{p\in\N}\,\pb_*\OP(p).$$
\sk
We will show that (3.2.1) is an isomorphism in Proposition~(3.5) below by using the Koszul and {\Cech} complexes as follows.
\ms\nin
{\bf 3.3.~Koszul and {\Cech} complexes.} Taking a local trivialization of the vector bundle $V$, and shrinking $X$, assume $V$ trivial.
Let $\xi_1,\dots,\xi_n$ be free generators of $\E^{\vee}$.
Set
$$\J_k:=(\xi_1^k,\dots,\xi_n^k)\subset\AX=\OX[\xi_1,\dots,\xi_n].$$
Then $\I^k$ in (3.1.1--2) can be replaced by $\J_k$.
Using (1.4.2), we get a free resolution of $\J_k$ by a complex whose $i$-th component is
$$\mopl_{|I|=1-i}\,\AX\xi_I^k,$$
where the $I$ are subsets of $\{1,\dots,n\}$, and $\xi_I:=\mprod_{i\in I}\xi_i$.
The differential is induced by the natural inclusions up to appropriate signs.
\sk
The derived graded localization $\R\Gamma(*\I)(\Gs)_{\ssb}$ for a quasi-coherent graded $\AX$-module $\Gs$ is represented by the inductive limit of the complex $K_k^{\prime\,\ssb}(\Gs)$ whose $i$-th component is
$$K_k^{\prime\,i}(\Gs)=\mopl_{|I|=i+1}\,\Gs(k|I|)\, \la\xi_I^{-k}\ra,
\leqno(3.3.1)$$
where $(k|I|)$ denotes the shift of degree as in (1.2.4), and the $\la\xi_I^{-k}\ra$ are formal symbols which are used to define the differential of $K_k^{\prime\,\ssb}(\Gs)$.
Indeed, the latter is induced by the morphisms
$$\xi_i^k:\Gs(k|I|)\,\la\xi_I^{-k}\ra\to\Gs(k|I'|)\, \la\xi_{I'}^{-k}\ra\q\h{for $\,I'=I\coprod\{i\}$},
\leqno(3.3.2)$$
with appropriate signs like usual {\Cech} complexes.
The formal symbols are also useful to define the canonical morphism between the $K_k^{\prime\,\ssb}(\Gs)$, which is induced by
$$\xi_i^{k'-k}:\Gs(k|I|)\,\la\xi_I^{-k}\ra\to \Gs(k'|I|)\, \la\xi_I^{-k'}\ra.
\leqno(3.3.3)$$
We then get
$${\rm ind\,lim}_{k\to\infty}\,\Gs(k|I|)\,\la\xi_I^{-k}\ra= \Gs[\xi_I^{-1}],$$
where the right-hand side is the graded localization of the graded $\AX$-module $\Gs$ by $\xi_I$.
\sk
As a conclusion, we have the following.
\ms\nin
{\bf 3.4.~Proposition.} {\it The derived graded localization $\R\Gamma(*\I)(\Gs)_{\ssb}$ for a quasi-coherent graded $\AX$-module $\Gs$ is represented locally on $X$ by the algebraic {\Cech} complex $K^{\prime\,\ssb}(\Gs)$ whose $i$-th component is}
$$K^{\prime\,i}(\Gs)=\mopl_{|I|=i+1}\,\Gs[\xi_I^{-1}].
\leqno(3.4.1)$$
\ms
As a corollary, we get
\ms\nin
{\bf 3.5.~Proposition.} {\it With the notation of {\rm (3.2)}, the morphism {\rm (3.2.1)} is an isomorphism for any quasi-coherent graded $\AX$-module $\Gs$, i.e.}
$$\Gs(*\I)_{\ssb}\simto\pb_*\Gt(\mb):=\mopl_{p\in\Z}\,\pb_*\Gt(p).$$
\ms\nin
{\it Proof.} Since the assertion is local we may assume $V$ trivial.
Then the assertion follows from Proposition~(3.4).
\ms\nin
{\bf 3.6.~Example.} Consider the simplest case where $X={\rm Spec}\,k$ with $k$ a field, $n=2$, and $\Gs=\AX$ so that $\PP=\PPP^1$, $\Gt=\OP$ and $\A=k[x,y]$.
Then the graded {\Cech} complex is given by
$$k[x,x^{-1},y]\oplus k[x,y,y^{-1}]\to k[x,x^{-1},y,y^{-1}],$$
and its graded cohomology groups at complex degree 0 and 1 are respectively
$$k[x,y],\q k[x^{-1},y^{-1}]\,\h{$\frac{1}{xy}$}.$$
This is compatible with a well-known calculation of $H^{\ssb}(\PPP^1,\Oc_{\PPP^1}(i))\,\,(i\in\Z)$.
\ms
Let $D^bG(\AX)_{qc}$ denote the derived category of bounded complexes of quasi-coherent graded $\AX$-modules.
\ms\nin
{\bf 3.7.~Theorem.} {\it For $\Gs\in D^bG(\AX)_{qc}$, let $\Gt$ denote the associated sheaf complex on $\PP$.
Then there is a canonical isomorphism}
$$\R\Gamma(*\I)(\Gs)_{\ssb}\cong\R\pb_*\Gt(\mb):= \mopl_{p\in\Z}\,\R\pb_*\Gt(p)\q\h{in}\,\,\,D^bG(\AX)_{qc}.$$
\ms\nin
{\it Proof.} Taking an injective resolution, we may assume that $\Gs$ is represented by a complex bounded below whose components are injective quasi-coherent graded $\AX$-modules.
Then the associated sheaf complex $\Gt$ is a complex of injective quasi-coherent $\OP$-modules.
So the assertion follows from Proposition~(3.5).
\sk
Combining Theorem~(3.7) with the distinguished triangle (3.1.3), we get the following.
\ms\nin
{\bf 3.8.~Theorem.} {\it For $\Gs\in D^bG(\AX)_{qc}$, there is a canonical distinguished triangle}
$$\R\Gamma_{\I}\1(\Gs)_{\ssb}\to\Gs\to\R\pb_*\Gt(\mb) \buildrel{+1}\over\to\q\h{in}\,\,\,D^bG(\AX).
\leqno(3.8.1)$$
\ms
It is also possible to prove directly the long exact sequence associated to (3.8.1) in a special case as below.
This is somehow related to \cite[Proposition~2.1]{Sch1} as mentioned in Remark~(3.11)(iii) below.
\ms\nin
{\bf 3.9.~Proposition.} {\it For a bounded complex of coherent graded $\AX$-module $\Gs$, there is a canonical long exact sequence of graded $\AX$-modules
$$\cdots\to\Hc^i_{\I}\1(\Gs)_{\ssb}\to\Hc^i\Gs\to R^i\pb_*\Gt(\mb)\to\Hc^{i+1}_{\I}(\Gs)_{\ssb}\to\cdots,
\leqno(3.9.1)$$
where the local cohomology sheaves $\Hc^i_{\I}\1(\Gs)_{\ssb}$ are defined by using a free resolution locally.
\sk
If $\Gs$ is a coherent graded $\AX$-module, i.e.\ if $\Hc^i\Gs=0$ for $i\ne 0$, then $(3.9.1)$ becomes isomorphisms
$$R^{i-1}\pb_*\Gt(p)\simto\Hc^i_{\I}\1(\Gs)_p\q(p\in\Z,\,i\ges 2)
\leqno(3.9.2)$$
and long exact sequences
$$0\to\Gamma_{\I}\1(\Gs)_p\buildrel{\alpha_p}\over\too\G_p \buildrel{\beta_p}\over\too\pb_*\Gt(p)\buildrel{\gamma_p}\over\too \Hc^1_{\I}\1(\Gs)_p\to 0\q(p\in\Z),
\leqno(3.9.3)$$
where $\alpha_p=0$, $\gamma_p=0$ and $\beta_p$ is bijective if $p\gg 0$, and $\alpha_p=0$, $\beta_p=0$ and $\gamma_p$ is bijective if $p\ll 0$.}
\ms\nin
{\it Proof.} It is sufficient to show the long exact sequence (3.9.1).
(Indeed, $\Hc^i_{\I}\1(\Gs)_p=0$ for $p\gg 0$ by Theorem~2 or using the arguments below.) Since the assertion is essentially local, we may assume that there is a quasi-isomorphism
$$\Ls\to\Gs,$$
where $\Ls^i=\A_X\motim_{\OX}\Ls^{\prime\,i}\,\,\,(i\in\Z_{\les 0})$ with $\Ls^{\prime\,i}$ a free graded $\OX$-module of finite rank.
The local cohomology sheaves can be defined by
$$\Hc^i_{\I}\,(\Gs)_{\ssb}:=\Hc^{i-n}(\Hc^n_{\I}\,(\Ls)_{\ssb}),$$
since $\Hc^k_{\I}\,(\Ls^i)_{\ssb}=0$ for $k\ne n$.
In the notation of (3.2), we have moreover
$$\Hc^n_{\I}\,(\Ls^i)\sssb=R^{n-1}\pb_*\Lt^i(\mb),$$
and there is a spectral sequence of graded $\AX$-modules
$$E_1^{p,q}=R^q\pb_*\Lt^p(\mb)\Longrightarrow R^{p+q}\pb_*\Gt(\mb).
\leqno(3.9.4)$$
Here $E_1^{p,q}=0$ unless $q=0$ or $n-1$, and we have
$$E_1^{p,0}=\pb_*\Lt^p(\mb)=\Ls^p,\q E_1^{p,n-1}=R^{n-1}\pb_*\Lt^p(\mb)=\Hc^n_{\I}\,(\Ls^p)\sssb.$$
Then
$$E_2^{p,0}=\Hc^p\Gs,\q E_2^{p,n-1}=\Hc^{p+n}_{\I}\,(\Gs)_{\ssb}.
\leqno(3.9.5)$$
We have $d_r^{p,q}=0$ except for
$$\h{$r=1\,$ with $\,q=0\,$ or $\,n-1,\,\,\,$ or $\,\,\,r=n\,$ with $\,q=n-1$.}$$
So the spectral sequence induces a long exact sequence
$$\cdots\to E_2^{p-n,n-1}\to E_2^{p,0}\to R^p\pb_*\Gt(\mb)\to E_2^{p+1-n,n-1}\to E_2^{p+1,0}\to\cdots.$$
This gives (3.9.1) by (3.9.5).
It is independent of the choice of the resolution, and is hence globally well-defined.
So Proposition~(3.9) follows.
\ms
As a corollary of Theorem~(3.8) or Proposition~(3.9) together with Theorem~2 in the introduction, we get the following.
\ms\nin
{\bf 3.10.~Corollary} (\cite{Sch1}).
{\it Let $(\M,F)$ be a Cohen-Macaulay holonomic filtered $\DX$-module.
Let $\PP$ be the projective cotangent bundle $P^*X:=(T^*X\setminus X)/\CC^*$, and $Z$ be the characteristic variety of $\M$ in $\PP$.
Let $d_x:=\dim Z\cap\PP_x$ for $x\in X$.
Then}
$$\cExt_{\OX}^j(\Gr^F_p\M,\OX)_x=0\q\h{for any}\,\,\,p\in\Z,\,\,\, \h{if}\,\,\,j>d_x+1.
\leqno(3.10.1)$$
\ms\nin
{\it Proof.} This follows from Theorem~2 and Theorem~(3.8) or Proposition~(3.9) applied to $\Gs:=\Gr^F_{\ssb}\M'$ with $(\M',F):=\DD(\M,F)$.
\ms\nin
{\bf 3.11.~Remarks.} (i) Proposition~(3.4) together with the definition of the algebraic {\Cech} complex $K^{\prime\,\ssb}(\Gs)$ in (3.3) can be generalized to the case $\Gs\in D^bG(\AX)_{qc}$ using the double complex construction.
\ms
(ii) It is well known that the category of quasi-coherent $\OP$-modules $M(\OP)_{qc}$ is equivalent to the category of quasi-coherent graded $\AX$-modules $G(\AX)_{qc}$ modulo the full subcategory $G_{\I}(\AX)_{qc}$ consisting of $\Gs'$ such that $\Gamma_{\I}\,(\Gs')\sssb\simto\Gs'$.
The inverse functor is given by the right-hand side of (3.2.1).
Proposition~(3.5) says that the functor $\Gamma(*\I)$ gives a canonical representative among $\Gs\in G(\AX)_{qc}$ corresponding to $\F\in M(\OP)_{qc}$, since $\R\Gamma(*\I)(\Gs')\sssb=0$ for $\Gs'\in G_{\I}(\AX)_{qc}$ by Proposition~(3.4).
\ms
(iii) When we consider quasi-coherent $\Oc_V$-modules corresponding to quasi-coherent $\AX$-modules, the information of the grading is lost.
For instance, the twisted sheaves $\OP(p)$ correspond to the shifted graded $\A$-modules $\A(p)$, but the pull-backs of $\OP(p)$ by the projection $\rho:V\setminus X\to\PP$ are independent of $p\in\Z$ where $X$ is identified with the zero-section.
(Indeed, in case $V=X\times\CC^n$, $\rho^*\OP(p)$ corresponds to a principal divisor on $X\times(\CC^n\setminus\{0\})$ defined by $x_1^p$ where $x_1$ is a coordinate of $\CC^n$.)
\sk
If we denote by $\Gt_V$ the coherent sheaf on $V$ associated to a coherent graded $\AX$-module $\Gs$, then, by the definition of $\Gt$ and $\Gt_V$, there are natural isomorphisms
$$\rho_*\rho^*\Gt\cong\mopl_{p\in\Z}\,\Gt(p),\q\rho^*\Gt\cong j^*\Gt_V,$$
where $j:V\setminus X\into V$ is the inclusion.
However, it is not easy to recover the information of the grading if we use the last isomorphism.
\sk
Proposition~2.1 in \cite{Sch1} was actually proved without the grading, and the grading on the local cohomology sheaves is essentially defined by using the exact sequence and the isomorphisms there.
In this case, it is not completely trivial to show the coincidence of this grading with the one obtained by taking a locally free resolution $\Ls\to\Gs$.
This can be shown, for instance, by using Proposition~(3.9).
Note also that the argument in loc.~cit.\ is sufficient for the proof of Corollary~(3.10) since we can take the direct sum over $p\in\Z$.
\ms
(iv) It is not quite trivial to show that the isomorphism in \cite{Sch1}, Lemma~2.6 is completely canonical.
In the notation of (3.2) in our paper, this can be proved, for instance, by applying duality to
$$R^{n-1}\pb_*\Lt(\mb):=\mopl_{p\in\Z}\,R^{n-1}\pb_*\Lt(p),$$
where $\Lt$ is the sheaf on $\PP$ associated to a locally free graded $\AX$-module $\Ls=\AX\sotim\Lc'$ with $\Lc'$ a locally free $\OX$-module.
In fact, $R^{n-1}\pb_*\Lt(\mb)$ is essentially used as a replacement for the graded local cohomology sheaf of $\Ls$ in loc.~cit.
The compatibility with the action of $\AX$ follows from the fact that the duality isomorphism is induced by the canonical pairing
$$\pb_*\Lt^{\vee}(\mb)\sotim_{\OX}R^{n-1}\pb_*\Lt(\mb)\to R^{n-1}\pb_*\Oc_{\PP}(\mb),$$
where $\Lt^{\vee}:=\Hcom_{\Oc}(\Lt,\Oc)$.
Here it is rather difficult to show the well-definedness of the isomorphism by using only the {\Cech} calculation since the {\Cech} covering itself may depend on the choice of the local trivialization of the cotangent bundle.
\ms
(v) Admitting the above arguments, the main theorem in loc.~cit.\ was essentially proved there not only for holonomic filtered $\D$-modules underling mixed Hodge modules, but also for any self-dual Cohen-Macaulay holonomic filtered $\D$-modules on smooth complex algebraic varieties.
\ms
(vi) If a Cohen-Macaulay holonomic filtered $\DX$-module $(\M,F)$ is self-dual up to the shift of the filtration by $w$, then the associated $\OP$-module $\Gt$ is Cohen-Macaulay and self-dual up to the twist by $\OP(-w)$.
If furthermore $Z:={\rm supp}\,\Gt$ (with reduced structure) is smooth and the multiplicity of $\Gt$ at the generic points is 1, then $\Gt$ is a locally free $\Oc_Z$-module of rank 1.
This can be shown for example by taking a local generic projection $U\to Z'$ inducing a finite \'etale morphism $U\cap Z\to Z'$ where $U\subset X$ is an open subvariety and $Z'$ is smooth.
\ms
(vii) If a Cohen-Macaulay filtered $\DX$-module $(\M,F)$ is supported on a strictly smaller subvariety $Z$ of $X$, then $\Gamma_{\I}(\Gr^F_{\ssb}\M)\sssb=0$.
Indeed, it corresponds to a coherent subsheaf of the coherent sheaf corresponding to $\Gr^F_{\ssb}\M$ whose support has dimension $\ge\dim X$ by the involutivity, and the subsheaf is supported on $Z\subset X\subset T^*X$.
So it vanishes by the definition of Cohen-Macaulay modules using depth.
\bs\bs
\centerline{\bf 4.~Examples}
\bs\nin
In this section we calculate some explicit examples to understand Theorem~2 and Proposition~(3.9) more concretely.
\ms\nin
{\bf 4.1.~Example} (Case $n=1$).
Let $X$ be an open disk $\Delta$ with coordinate $x$ where $n=1$.
For simplicity, assume $(\M,F)|_{\Delta^*}$ is a nilpotent orbit of weight $w$ with unipotent monodromy and one Jordan block of size $w+1$.
Let $e_0,\dots,e_w$ be a basis of the Deligne extension such that $e_j=(x\rd_x)^{w-j}e_w$ where $\rd_x:=\rd/\rd x$.
Assume the Hodge filtration $F^p$ on $\M|_{\Delta^*}$ is generated by $e_j\,(j\ges p)$, and $(\M,F)$ is the filtered left $\DX$-module underlying a pure Hodge module with strict support $X$.
The latter has weight $w+1$, and hence
$$(\M',F):=\DD(\M,F)\cong(\M,F)(w+1),$$
see (2.2) for $\DD(\M,F)$.
Set $v_p=\rd_x^{w-p} e_w\in F_{-p}\M$.
Let $[v_p]\in\Gr^F_{-p}\,\M$ denote the class of $v_p$ in the graded quotient, and similarly for $[e_p]$, etc.
Then
$$\Gr^F_{w+1-p}\,\M'\cong\Gr^F_{-p}\,\M=\begin{cases}\,0&\h{if}\,\,\,\,p>w,\\ \OX[v_p]\cong\OX&\h{if}\,\,\,\,p\in[0,w],\\ \OX[v_p]\cong\OX/(x^w)&\h{if}\,\,\,\,p<0,\end{cases}
\leqno(4.1.1)$$
Set $\xi=\Gr^F_1\rd_x\in\Gr^F_1\DX$.
By (4.1.1) we get for $k\gg 0$
$$\Ker(\xi^k:\Gr^F_{-p}\,\M\to\Gr^F_{k-p}\,\M)=\begin{cases}\OX[x^wv_p]& \h{if}\,\,\,\,p\in[0,w],\\ \,0&\h{if}\,\,\,\,p\notin[0,w].\end{cases}
\leqno(4.1.2)$$
These are compatible with Theorem~1 since
$$[v_p]=[x^{p-w}e_p]\in\Gr^F_{-p}\,\M|_{\Delta^*},\q [x^wv_{w-p}]=[x^{w-p}e_{w-p}]\in\Gr^F_{p-w}\,\M|_{\Delta^*}.$$
Note that there is a self-pairing $\la *,*\ra$ of $\M|_{\Delta^*}$ corresponding to the above self-duality, and
$$\la e_i,e_j\ra=(-1)^i\delta_{i+j-w}.$$
Indeed, we have $e_j=\exp\bl(-(\log x)N\br)u_j$ with $u_j$ locally constant multivalued sections and $N:=-\log T/2\pi i$.
Then $u_j=N^{w-j}u_w$, and $\la e_i,e_j\ra=\la u_i,u_j\ra$ since $\la Nu,v\ra=-\la u,Nv\ra$.
\sk
Set $\Gs:=\Gr^F_{\ssb}\M$.
Since $\PPP^0=pt$, we have for any $p\in\Z$
$$R^i\pb_*\Gt(p)=\begin{cases}\OX/(x^w)&\h{if}\,\,\,\,i=0,\\ \,0&\h{if}\,\,\,\,i\ne 0.\end{cases}$$
Set $R:=\Oc_{X,0}$, $\af:=(x^w)\subset R$, $Q:=R/\af$.
Then (3.9.3) is identified with
$$\aligned 0\to 0\to Q\to Q\to 0\to 0&\q\h{if}\,\,\,p>0,\\ 0\to\af\to R\to Q\to 0\to 0&\q\h{if}\,\,\,p\in[-w,0],\,\\ 0\to 0\to 0\to Q\to Q\to 0&\q\h{if}\,\,\,p<-w.\endaligned
\leqno(4.1.3)$$
Note that $(\Hc^1_{\I}\,\Gs)_p$ is identified with $\Ext_{\OX}^1(\Gr^F_{-p-w}\M,\OX)$ by Theorem~2 and that this is compatible with the isomorphism
$$\Ext^1_R(Q,R)\cong Q.$$
\ms\nin
{\bf 4.2.~Remark.} Let $(\M,F)$ be a filtered holonomic $\DX$-module underlying a polarizable Hodge module with strict support $X$ where $X$ is a complex manifold.
Let $D$ be a {\it smooth} hypersurface of $X$ defined by a coordinate function $x_1$.
Set $U=X\setminus D$.
Assume $(\M,F)|_U$ is a variation of Hodge structure of odd weight $n$.
Set $m=[n/2]$ so that $n=2m+1$.
Let $T$ be the local monodromy around $D$.
Assume $T$ is unipotent and ${\rm rank}\,(T-id)=1$.
Then the reader might expect that the same argument as above would imply locally on $X$
$$\cExt_{\OX}^1(\Gr^F_p\M,\OX)\cong\Oc_D\q\h{for}\,\,\,p>-m,
\leqno(4.2.1)$$
contradicting the reflexivity in the universal family case, see \cite{Sch2}.
However, the same argument does not seem to hold in the universal family case as is seen below:
\sk
Consider the exact sequence (3.9.3).
Its last term corresponds to the left-hand side of (4.2.1) by Theorem~2.
We have locally on $X$
$$\pb_*\Gt(p)\cong\Oc_D\q\h{for any}\,\,\,p\in\Z.$$
However, for $p\in[-n,-m-1]$, $\gamma_p$ in (3.9.3) may vanish, i.e., $\beta_p$ may be surjective, and the cokernel of $\alpha_p$ may be locally isomorphic to $\Oc_D$ in the universal family case.
(Indeed, let $(u_i)$ be a basis of multivalued horizontal sections on $U$ such that $Nu_i=0$ for $i\ne 1$ and $Nu_1=u_2$, shrinking $X$ if necessary.
Then $\beta_p$ is surjective if there is a nonzero element of $\G_p=\Gr^F_p\M$ represented by $\sum_i g_iu_i$ such that $\rd_{x_1}^kg_1|_D\ne 0$ with $p=k-m$, where $g_i\in\OX\,\,(i\ne 2)$ and $g_2\in\OX[\log x_1]$.)
\ms\nin
{\bf 4.3.~Example} (Supported on a hypersurface with an ordinary double point).
Assume $X=\AAA^n_k$ and $n\ges 4$, where $k$ is a subfield of $\CC$.
Let $x_1,\dots,x_n$ be the coordinates of $\AAA^n_k$. Set
$$f=\msum_{i=1}^n\,x_i^2,\q Z=\{f=0\}\subset X,$$
and
$$\Nc=\OX(*Z),\q\Ncb=\OX(*Z)/\OX.$$
These $\DX$-modules have the Hodge filtration $F$ so that they underly mixed Hodge modules.
Let $(\M,F)$ be the filtered $\DX$-module underlying the pure Hodge module ${\rm IC}_Z\Q(-1)$ of weight $n+1$ where $(-1)$ is the Tate twist.
Let $m:=[n/2]$.
It is well known (see for instance \cite{DSW}, \cite{hfil}, etc.) that there is a natural inclusion
$$\iota:(\M,F)\into(\Ncb,F),$$
underlying an injection of mixed Hodge modules.
Moreover,
$${\rm Coker}\,\iota=\begin{cases}(\B,F[m-1])&\h{if $\,n\,$ is even,}\\ \,0&\h{if $\,n\,$ is odd.}\end{cases}$$
Here $\B:= k[\rd_1,\dots,\rd_n]$ with $\rd_i:=\rd/\rd x_i$, and it has the filtration $F$ defined by the degree of polynomials in $\rd_1,\dots,\rd_n$.
\sk
Since $(\M,F)$ is self-dual up to a Tate twist, and has weight $n+1$, we have
$$\M':=\DD(\M,F)=(\M,F[n+1]),\,\,\h{i.e.}\q \Gr^F_{\ssb}\M'=\Gr^F_{\ssb-n-1}\M.$$
Then Theorem~2 implies the duality isomorphisms for $q\in\Z$
$$\R\Gamma_{\I}(\Gr^F_{\ssb}\M)_q=\R\Hcom_{\OX}(\Gr^F_{-1-q}\,\M,\OX).
\leqno(4.3.1)$$
\sk
We now calculate both sides of (4.3.1) by using (3.9.2--3) for the left-hand side.
Let
$$\aligned &M:=\Gamma(X,\M),\q N:=\Gamma(X,\Nc), \q R:=\Gamma(X,\OX),\q \Rb:=\Gamma(Z,\Oc_Z),\\ &\q\h{so that}\q N=R[f^{-1}],\q R=k[x_1,\dots,x_n],\q\Rb=R/fR.\endaligned$$
Set $R^{\ges p}:=\msum_{|\nu|\ges p}\,k\,x^{\nu}\subset R$ with $x^{\nu}:=\prod_{i=1}^nx_i^{\nu_i}$ for $\nu=(\nu_1,\dots,\nu_n)\in\N^n$.
(Similarly for $R^p$, $R^{\les p}$ with $|\nu|\ges p$ replaced by $|\nu|=p$ or $|\nu|\les p$.) Let $\Rb^{\ges p}$, etc.\ denote the image of $R^{\ges p}$ in $\Rb$.
By \cite[Lemma~1.5]{DSW}, we have
$$F_p(R[f^{-1}])=\begin{cases}f^{-p-1}R^{\ges p-m+1}\subset R[f^{-1}]& \h{if}\,\,\,\,p\ges 0,\\ \,0&\h{if}\,\,\,\,p<0.\end{cases}$$
After some calculation we then get the short exact sequences for $p\ges 0:$
$$\aligned 0\to R^{\ges p-m}\buildrel{f}\over\to R^{\ges p-m+1}\to \Gr^F_pM\to 0&\q\h{if $\,n\,$ is odd},\\ 0\to R^{\ges p-m}\buildrel{f}\over\to R^{\ges p-m+2}\to \Gr^F_pM\to 0&\q\h{if $\,n\,$ is even.}\endaligned
\leqno(4.3.2)$$
\sk
We first consider the case $n$ odd, i.e., $n=2m+1$.
For $p\ges 0$, set
$$\DD^j_R\Gr^F_pM:=\Ext_R^j(\Gr^F_pM,R).$$
This vanishes for $p<0$.
Calculating the dual by using (4.3.2), we get for $p\ges 0$
$$\DD^j_R\Gr^F_pM=\begin{cases}\,\Rb&\h{if}\,\,\,\,j=1,\\ \bl(\Rb^{\,\les p-m}\br)^{\vee}&\h{if}\,\,\,\,j=n-1,\\ \bl(R^{\,p-m-1}\br)^{\vee}&\h{if}\,\,\,\,j=n,\\ \,0&\h{otherwise}.\end{cases}
\leqno(4.3.3)$$
\sk
Let $G^{\ssb}$ be the $\mf$-adic filtration on $\Gs:=\Gr^F_{\ssb}\M$ where $\mf$ is the maximal ideal of $R$ generated by $x_1,\dots,x_n$.
The induced filtration on $\Gt$ will be denoted also by $G^{\ssb}$.
Identifying $R^j\pb_*\Gr_G^k\,\Gt(q)$ with $H^j(\PPP^{n-1},\Gr_G^k\,\Gt(q))$, we then get for $q\in\Z$
$$R^j\pb_*\Gr_G^k\,\Gt(q)=\begin{cases}R^{\,q-m+1}&\h{if}\,\,\,\,j=0,\,\,k=0, \\ \Rb^{\,q-m+1+k}&\h{if}\,\,\,\,j=0,\,\,k\ges 1,\\ \bl(\Rb^{\,-q-m-k}\br)^{\vee}&\h{if}\,\,\,\,j=n-2,\,\,k\ges 1,\\ \bl(R^{\,-q-m-2}\br)^{\vee}&\h{if}\,\,\,\,j=n-1,\,\,k=0.\\ \,0&\h{otherwise}.\end{cases}
\leqno(4.3.4)$$
Here the differentials $d_r^{p,q}$ of the associated spectral sequence vanish (by using $n\ges 4$).
\sk
Set $p+q=-1$ in order to apply (4.3.1).
Then (4.3.4) agrees with (4.3.3) via (4.3.1) and (3.9.2--3).
Indeed, $\Gamma_{\I}\,(\Gs)_{\ssb}=0$ by Remark~(3.11)(vii), and via the exact sequence (3.9.3), $\pb_*\Gt(q)$ contributes to $\G_q$ for $q\ges 0$, and to $\Hc^1_{\I}\1(\Gs)_q$ for $q\le-1$ (i.e. for $p\ges 0$).
\sk
In case $n$ is even, we get for $p\ges 0$
$$\DD^j_R\Gr^F_pM=\begin{cases}\Rb&\h{if}\,\,\,\,j=1,\\ \bl(\Rb^{\,\les p-m+1}\br)^{\vee}&\h{if}\,\,\,\,j=n-1,\\ \,0&\h{if}\,\,\,\,j\ne 1,\,n-1,\end{cases}
\leqno(4.3.5)$$
by using (4.3.2).
We have for $q\in\Z$
$$R^j\pb_*\Gr_G^k\,\Gt(q)=\begin{cases}\Rb^{\,q-m+2+k} &\h{if}\,\,\,\,j=0,\,\,k\ges 0,\\ \bl(\Rb^{\,-q-m-k}\br)^{\vee}&\h{if}\,\,\,\,j=n-2,\,\,k\ges 0,\\ \,0&\h{otherwise.}\end{cases}
\leqno(4.3.6)$$
Setting $p+q=-1$, we see similarly that (4.3.6) agrees with (4.3.5) via (4.3.1) and (3.9.2--3) in case $n$ is even.
(It is quite difficult to generalize the above calculation to the case of singularities of more complicated type in \cite[Theorem~0.7]{hfil}.)
\bs\bs
\centerline{\bf 5.~Vanishing of higher extension groups}
\bs\nin
In this section we explain an application of Corollary~(3.10) in a more general situation than in \cite[Section 1.3]{Sch1}.
\ms\nin
{\bf 5.1.~Notation.} The application to the hypersurface case in \cite[Section 1.3]{Sch1} is valid not only for universal families, but also for families of hypersurfaces whose total space is smooth.
More precisely, let $X$ be a smooth complex projective variety, and assume there is a family of hypersurfaces
$$\h{$\Y=\coprod_{s\in S}Y_s\into \X:=X\times S$},$$
such that the total space $\Y$ and the parameter space $S$ are smooth.
Let $f:\Y\to S$ be the composition with the projection $pr:\X\to S$.
Let $(\M,F)$ be a direct factor of the cohomological direct image $H^if_*(\Oc_{\Y},F)$ as a filtered left $\D$-module.
(The direct image is defined as in \cite{mhp} via the transformation in (2.1.3).) Let $\PP$ be the projective cotangent bundle $P^*S:=(T^*S\setminus S)/\CC^*$ where $S$ is identified with the zero-section of $T^*S$.
There is a canonical morphism
$$\phi_s:\Sing Y_s\to \PP_s,$$
such that
$$(P^*_{\Y}\X)_y=\{pr^*\phi_s(y)\}\,\,\,\h{in}\,\,\,P^*_y\X \q\h{for}\,\,\,y\in\Sing Y_s,
\leqno(5.1.1)$$
where $pr^*:\PP_s=P^*_sS\to P^*_y\X$ is induced by the projection $pr:\X\to S$, and $P^*_{\Y}\X\subset P^*\X$ is the projective conormal bundle of $\Y$ in $\X$.
We have the following (see \cite{Sch1}):
\ms\nin
{\bf 5.2.~Corollary.} {\it Let $(\M,F)$ be as above.
Let $s\in S$, $j\in\Z$ with $j>\dim\Sing Y_s+1$, or more precisely $j>\dim\phi_s(\Sing Y_s)+1$.
Then we have for any $p\in\Z$
$$\cExt_{\Oc_S}^j(\Gr^F_p\M,\Oc_S)_s=0,
\leqno(5.2.1)$$
where the left-hand side is the stalk at $s$ of the coherent sheaf $\cExt_{\Oc_S}^j(\Gr^F_p\M,\Oc_S)$.}
\ms\nin
{\it Proof.} This follows from Corollary~(3.10) (see also \cite{Sch1}) by using a well-known estimate of the characteristic variety of the direct image of a $\D$-module, which is essentially expressed by the above morphism $\phi_s$ in this special case.
(Here the assertion is independent of $i$ since the assumption about the direct factor is used only for the estimate of the characteristic variety of $\M$.) This finishes the proof of Corollary~(5.2).
\ms\nin
{\bf 5.3.~Remark.} In examples, we usually have the equality
$$\dim\phi_s(\Sing Y_s)=\dim\Sing Y_s,
\leqno(5.3.1)$$
as in Example~(5.4) below.
It is not easy to get an example where the equality does not hold, except for the case where the family comes from hypersurfaces of a lower dimensional variety $X'$ via a smooth morphism $X\to X'$.
\ms\nin
{\bf 5.4.~Example.} Let $X,S$ be projective space $\PPP^n$ with projective coordinates respectively $x_0,\dots,x_n$ and $s_0,\dots,s_n$.
Consider a family of hypersurfaces defined by
$$Y_s:=\{\msum_{i=0}^n\,s_ix_i^2=0\}\subset X.$$
For $I\subset \{0,\dots,n\}$ with $I\ne\emptyset$, set $I^c:=\{0,\dots,n\}\setminus I$, and
$$S_I:=\{\h{$s_i\ne 0\,$ if $\,i\in I,\,$ and $\,s_i=0\,$ if $\,i\in I^c$}\}\subset S.$$
Then
$$\Sing Y_s=X^I:=\{x_i=0\,(i\in I)\}\subset X\q\h{for}\,\,\, s\in S_I,$$
since we have on $(\CC^{n+1}\setminus\{0\})^2$
$$d(\msum_{i=0}^n\,s_ix_i^2)=\msum_{i=0}^n\,2s_ix_idx_i +\msum_{i=0}^n\,x_i^2ds_i.$$
The latter also implies that the morphism $\phi_s$ in (5.1.1) is expressed for $s\in S_I$ by
$$\Sing Y_s=X^I\ni[x_i]_{i\in I^c}\mapsto \bl[x_i^2\br]_{i\in I^c}\in\PPP^{n-|I|},
\leqno(5.4.1)$$
where we use a trivialization of the projective conormal bundle
$$P^*_{S_I}S\cong\PPP^{n-|I|}\times S_I,$$
given by the sections $d(s_i/s_{i_0})\,(i\in I^c)$ choosing some $i_0\in I$.
\ms
We have $\dim S_I=|I|-1$, $\dim X^I=n-|I|$, and for $s\in S_I$, $j\in[0,2n-2]$
$$H^j(Y_s,\Q)\cong\begin{cases}\buildrel{2}\over\mopl\Q(-j/2)& \h{if $\,j=2n-|I|\,$ with $\,j\,$ even,}\\ \,\Q(-j/2)&\h{if $\,j\ne 2n-|I|\,$ with $\,j\,$ even,}\\ \,\,0&\h{if $\,j\,$ is odd.}\end{cases}
\leqno(5.4.2)$$
This can be shown for example by increasing induction on $n\ges |I|-1$ using the projection
$$Y_s\setminus P_{i_0}\to\Ybs:=\{\msum_{i\ne i_0}s_ix_i^2=0\} \subset\PPP^{n-1},$$
where $P_{i_0}:=\{x_i=0\,(i\ne i_0)\}\subset Y_s$ with $i_0\in I^c$.
Indeed, we have an isomorphism
$$H^j_c(Y_s\setminus P_{i_0},\Q)\simto\Hct^j(Y_s,\Q),$$
where $\Hct^j$ is the reduced cohomology, and $Y_s\setminus P_{i_0}$ is a line bundle over $\Ybs$ so that
$$H^j_c(Y_s\setminus P_{i_0},\Q)=H^{j-2}(\Ybs,\Q)(-1).$$
\sk
The above calculation implies the decomposition
$$\R f_*\Q_{\Y}[2n-1]\cong \bl(\mopl_{i=0}^{n-1}\,(\Q_S[n])(-i)[n-1-2i]\br)\oplus \bl(\mopl_{|I|:\,{\rm even}}\,(j_I)_!L_I[d_I]\br),
\leqno(5.4.3)$$
where $j_I:S_I\into S$ is the canonical inclusion, $L_I$ is a $\Q$-local system of rank 1 on $S_I$, and $d_I:=\dim S_I=|I|-1$.
Moreover, the local monodromy of $L_I$ around $S_{J}$ with $|J|=|I|-1$ is given by multiplication by $-1$.
Hence
$$(j_I)_!L_I[d_I]=\R(j_I)_*L_I[d_I]=(j_I)_{!*}L_I[d_I],
\leqno(5.4.4)$$
where $(j_I)_{!*}$ is the intermediate direct image \cite{BBD}.
We then get
$$f_*(\Oc_{\Y},F)\cong \bl(\mopl_{i=0}^{n-1}\,(\Oc_S,F[-i])[n-1-2i]\br)\oplus \bl(\mopl_{|I|:\,{\rm even}}\,(\M_I,F)\br),
\leqno(5.4.5)$$
where the left-hand side is the direct image as a filtered left $\D$-module, and $(\M_I,F)$ is a filtered regular holonomic $\D_S$-module corresponding to $(j_I)_{!*}L_I[d_I]\sotim_{\Q}\CC$ by the de Rham functor $\DR_S$, and underlying a pure Hodge module of weight $2n-1$
$$\bl((\M_I,F),(j_I)_{!*}L_I[d_I]\br).$$
The Hodge filtration $F$ on $\M_I$ can be calculated locally on $S$ as in \cite[Proposition~3.5]{mhm}.
Locally we have
$$\Gr^F_p\M_I\cong\mopl_{J\subset I}\,\E_{p,I,J},$$
where the $\E_{p,I,J}$ are locally free $\Oc_{\Sb_J}$-modules with $\Sb_J:=\bigcup_{K\subset J}S_K$ (the closure of $S_J$), and $\E_{p,I,J}\ne 0$ for any $J\subset I$ if $|I|$ is even and $p\gg 0$.
Since $\Sb_J$ is smooth with pure codimension $n+1-|J|$, we have
$$\cExt_S^j(\E_{p,I,J},\Oc_S)\ne0\iff j=n+1-|J|,$$
and $\dim \Sing Y_s=n-|J|$ for $s\in S_J$.
So the estimate in Corollary~(5.2) is optimal.

\end{document}